\documentclass[11 pt]{article}

\usepackage{amssymb}

\def\nmark{\mbox{$\rm\bf\kern0.2em\rule{0.06em}{1.45ex}\kern-0.3em
N$}}
\def\dmark{\mbox{$\rm\bf\kern0.2em\rule{0.06em}{1.45ex}\kern-0.3em
D$}}
\def\cmark{\mbox{$\rm\bf\kern0.2em\rule{0.06em}{1.45ex}\kern-0.3em
C$}}
\def\rmark{\mbox{$\rm\bf\kern0.2em\rule{0.06em}{1.45ex}\kern-0.3em
R$}}

\date{30 January 2018}

\begin{document}
\title{\large \bf
Complex symmetric weighted composition operators}
 \author{ Mahsa Fatehi}

{\maketitle}
\begin{abstract}
In this paper we find all complex symmetric weighted composition operators with  special conjugations. Then we give spectral properties of these complex symmetric weighted composition operators.
\end{abstract}

\footnote{
{\it key words and phrases}:  Complex symmetric operator, conjugation, weighted composition operator.}

\section{Introduction}

Let $\mathbb{D}$ denote the open unit disk in the complex plane. The Hardy space, denoted $H^{2}(\mathbb{D})=H^{2}$, is  the set of all analytic functions $f$ on $\mathbb{D}$, satisfying the norm condition
 $$\|f\|^{2}=\lim _{r \rightarrow 1}\int_{0}^{2\pi}|f(re^{i\theta})|^{2}\frac{d\theta}{2\pi}< \infty.$$
 The space $H^{\infty}(\mathbb{D})=H^{\infty}$ consists of all the functions that are analytic and bounded on $\mathbb{D}$, with supremum norm $\|f\|_{\infty}=\sup_{z \in \mathbb{D}}|f(z)|$.\par

 Let $\varphi$ be an analytic map from the open unit disk $\mathbb{D}$ into itself. The operator that takes the analytic map $f$ to $f \circ \varphi$ is a composition operator and is denoted by $C_{\varphi}$. A natural generalization of a composition operator is an operator that takes $f$ to $\psi \cdot f \circ \varphi$, where $\psi$ is a fixed analytic map on $\mathbb{D}$. This operator is aptly named a weighted composition operator and is usually denoted by $C_{\psi,\varphi}$. More precisely, if $z$ is in the unit disk then
$(C_{\psi,\varphi}f)(z)=\psi(z)f(\varphi(z)).$

The automorphisms of $\mathbb{D}$, that is, the one-to-one analytic maps of the disk onto itself, are just the functions
$\varphi(z)=\lambda\frac{a-z}{1-\overline{a}z},$
where $|\lambda|=1$ and $|a| < 1$.

Let $P$ denote the orthogonal
projection of $L^{2}(\partial \mathbb{D})$ onto $H^{2}$. For each $b \in {L^{\infty}(\partial \mathbb{D})}$,
the Toeplitz  operator $T_{b}$ acts on $H^{2}$ by
$T_{b}(f)=P(bf)$.\par
In \cite{co2}, Cowen obtained an adjoint formula of a composition operator whose symbol is a linear-fractional self-map of $\mathbb{D}$; for $\varphi(z)=\frac{az+b}{cz+d}$ which is a linear-fractional self-map of $\mathbb{D}$, he showed that $C_{\varphi}^{\ast}=T_{g}C_{\sigma_{\varphi}}T_{h}^{\ast}$, where $\sigma_{\varphi},g$ and $h$ are the Cowen auxiliary functions given by $\sigma_{\varphi}(z):=({\overline{a}z-\overline{c}})/({-\overline{b}z+\overline{d}})$, $g(z):=(-\overline{b}z+\overline{d})^{-1}$ and
$h(z):=(cz+d)$.  We can see that $\sigma_{\varphi}(z)=\frac{1}{\overline{\varphi^{-1}(\frac{1}{\overline{z}})}}$, so $\sigma_{\varphi}$ maps $\mathbb{D}$ into itself. It follows that $\varphi(a)=a$ if and only if $\sigma_{\varphi}(\frac{1}{\overline{a}})=\frac{1}{\overline{a}}$, where $a \in \mathbb{C}$. Note that  $g$ and $h$ are in  $H^{\infty}$. If $\varphi(\zeta)=\eta$
for $\zeta,\eta \in \partial \mathbb{D}$, then
$\sigma_{\varphi}(\eta)=\zeta$. We know that $\varphi$ is an
automorphism if and only if $\sigma_{\varphi}$ is, and in this case
$\sigma_{\varphi}=\varphi^{-1}$. we
will write $\sigma$ for $\sigma_{\varphi}$ except when confusion
could arise. From now on, unless otherwise stated, we
assume that $\sigma$, $h$ and $g$ are given as above. \par

A bounded operator $T$ on a complex Hilbert space $H$ is said to be a complex symmetric operator if there exists a conjugation $C$ (an isometric, antilinear involution) such that $CT^{\ast}C=T$. The complex symmetric
operators class was initially addressed by Garcia and Putinar (see \cite{gar1} and \cite{gar2})
and includes the normal operators,
Hankel operators and Volterra integration operators. Invoking \cite[Theorem 2]{gar3}, any composition operator with an involutive automorphism symbol is complex symmetric. In  \cite{noor}, Bourdon et al. showed that among the automorphisms of $\mathbb{D}$, only the elliptic ones may introduce complex symmetric operators. Moreover, they proved that for $\varphi$, not the rotation and involutive automorphism, which is elliptic automorphism of order $q$ that $4\leq q\leq\infty$, $C_{\varphi}$ is not complex symmetric.  In this paper we use the symbol $J$ for the special conjugation that $(Jf)(z)=\overline{f(\overline{z})}$ for each analytic function $f$. In \cite{ham} and \cite{skl}, all $J$-symmetric weighted composition operators were characterized.  Recently in \cite{nst} Narayan et al. have found complex symmetric composition operators whose symbols are linear-fractional, but not an automorphism.\par
In the second section of this paper, first we find all unitary weighted composition operators which are $J$-symmetric. Then we consider the special conjugations which are the  products of these unitary weighted composition operators and the conjugation $J$. Next in Theorem 2.5, we obtain complex symmetric weighted composition operators with these conjugations. In addition, in Theorem 2.7, we characterize all complex symmetric weighted composition operators which are isometries.  \par
In the third section,  we provide a characterization of $\varphi$, when $\varphi$ has the form which was stated in \cite[Proposition 2.9]{ham} and \cite[Theorem 3.3]{skl}. Finally, we obtain the spectrum and spectral radius of some complex symmetic weighted composition operators that  study in the  second section.

\section{Weighted composition operators}

Suppose that $\varphi(z)=\frac{az+b}{1-cz}$ is a linear-fractional self-map of $\mathbb{D}$. If $\varphi$ is written as $\varphi(z)=a_{0}+\frac{a_{1}z}{1-a_{0}z}$, then it is not hard to see that $c=b=a_{0}$. We use this fact frequently in this paper.  \par
An operator $T$ is said to be unitary if $T^{\ast}T=TT^{\ast}=I$.
In the following proposition, we find all unitary weighted composition operators $C_{\psi,\varphi}$ which are  $J$-symmetric.  \\ \par

{\bf Proposition 2.1.} {\it The weighted composition operator  $C_{\psi,\varphi}$ is unitary and $J$-symmetric if and only if either $\psi(z)=c\frac{(1-|p|^{2})^{1/2}}{1-\overline{p}z}$ and $\varphi(z)=\frac{\overline{p}}{p}\frac{p-z}{1-\overline{p}z}$, where $p \in \mathbb{D}-\{0\}$ and $|c|=1$ or $\psi \equiv \mu$  and $\varphi(z)=\lambda z$, when $|\mu|=|\lambda|=1$.}\bigskip

{\bf Proof.} Let $C_{\psi,\varphi}$ be unitary and $J$-symmetric. By \cite[Theorem 6]{Bourdon},
\begin{equation}
\varphi(z)=\lambda \frac{p-z}{1-\overline{p}z}
\end{equation}
and
\begin{equation}
\psi(z)=c\frac{(1-|p|^{2})^{1/2}}{1-\overline{p}z},
\end{equation}
where $|\lambda|=|c|=1$. First suppose that $p \neq 0$. We can see that $\varphi(0)=\lambda p$ and $\varphi'(0)=\lambda (|p|^{2}-1)$. Let $\widetilde{\varphi}(z)=\lambda p+\frac{\lambda(|p|^{2}-1)z}{1-\lambda pz}$. It is not hard to see that $\widetilde{\varphi}\equiv \varphi$ if and only if $\lambda=\frac{\overline{p}}{p}$. Invoking \cite[Theorem 3.3]{skl} (see also \cite[Proposition 2.9]{ham}), the conclusion follows for $p\neq 0$. Letting $p=0$ in Equations (1) and (2), we get  $\psi$ is a constant function and $\varphi(z)=-\lambda z$, when $|\lambda|=1$.\par
Conversely, suppose that either $\varphi(z)=\frac{\overline{p}}{p}\frac{p-z}{1-\overline{p}z}$ and $\psi(z)=c\frac{(1-|p|^{2})^{1/2}}{1-\overline{p}z}$ or $\psi$ is a constant function and $\varphi(z)=\lambda z$. In this both cases, by \cite[Theorem 6]{Bourdon}, $C_{\psi,\varphi}$ is unitary. Rewriting $\varphi$, we see that $\varphi(z)=\frac{\overline{p}}{p}\frac{p-z}{1-\overline{p}z}=\overline{p}+\frac{\frac{\overline{p}}{p}(|p|^{2}-1)z}{1-\overline{p}z}.$
Again \cite[Proposition 2.9]{ham} and \cite[Theorem 3.3]{skl} imply that in these both cases the weighted composition operators $C_{\psi,\varphi}$ are $J$-symmetric.\hfill $\Box$ \\ \par

From now on, unless
otherwise stated, we assume that $\varphi_{p}(z)=\frac{\overline{p}}{p}\frac{p-z}{1-\overline{p}z}$, where $p \in \mathbb{D}-\{0\}$ and $\psi_{p}(z)=c\frac{(1-|p|^{2})^{1/2}}{1-\overline{p}z}$, when $p \in \mathbb{D}$ and $|c|=1$. The proof of the next lemma is left to the reader.\\ \par

{\bf Lemma 2.2.} {\it If $U$ is a unitary and complex symmetric operator with conjugation $C$, then $UC$ is a conjugation.}\bigskip

Let $C$ be a conjugation. Then $CJ=W$ is a unitary operator and $W$ is both $C$-symmetric and $J$-symmetric (see \cite[Lemma 3.2]{gar}). We have $C=WJ$. Then all conjugations can be considered as a product of a $J$-symmetric unitary  operator $W$ and the conjugation $J$.\\ \par

{\bf Proposition 2.3.} {\it Suppose that $U$ is unitary and complex symmetric with conjugation $WJ$, where $W$ is unitary. Then an operator $A$ is $WJ$-symmetric if and only if $UA$ is $UWJ$-symmetric. }\bigskip

{\bf Proof.} Suppose that $A$ is $WJ$-symmetric. Invoking Lemma 2.2,  $UWJ$ is a conjugation. We have
$$UWJ(UA)^{\ast}UWJ=UWJA^{\ast}U^{\ast}UWJ=UA,$$
so $UA$ is complex symmetric with conjugation $UWJ$.\par
Conversely, suppose that $UA$ is $UWJ$-symmetric. We see that
$$WJA^{\ast}WJ=U^{\ast}UWJ(UA)^{\ast}UWJ=U^{\ast}UA=A.$$
Hence, $A$ is $WJ$-symmetric.\hfill $\Box$ \\
\par

Assume that an operator $B$ is $UJ$-symmetric. Let $W=I$ in Proposition 2.3. By Proposition 2.3, $U^{\ast}B$ is $J$-symmetric. We have $B=UU^{\ast}B$.  It shows that every complex symmetric operator can be written as a product of a unitary $J$-symmetric operator and a $J$-symmetric operator. Since recently a lot of $J$-symmetric operators have been found, this idea may be useful in order to obtain complex symmetric operators more. \par

 In the following example, we find a complex symmetric Toeplitz operator $T_{f}$ with $|f|=1$ on $\partial \mathbb{D}$ (see \cite[Corollary 2.2]{kolee}). \\ \par

{\bf Example 2.4.} Suppose that $p\in (-1,1)$ is a real number. By \cite[Lemma 2.1]{noor} and  \cite[Lemma 2.2]{noor}, $U_{\varphi_{p}}$ is $J$-symmetric, when $U_{\varphi_{p}}$ is the unitary part in the polar decomposition of $C_{\varphi_{p}}$ (note that in this case $\varphi_{p}$ is an involutive automorphism). By the proof of \cite[Lemma 4.7]{mw2}, $U_{\varphi_{p}}=C_{\varphi_{p}}T_{\frac{|1-pz|}{(1-|p|^{2})^{1/2}}}$.  Proposition 2.1 implies that  $C_{\psi_{p},\varphi_{p}}$ is $J$-symmetric. Invoking Proposition 2.3, $C_{\psi_{p},\varphi_{p}}C_{\varphi_{p}}T_{\frac{|1-pz|}{(1-|p|^{2})^{1/2}}}$ is $C_{\psi_{p},\varphi_{p}}J$-symmetric. Then $T_{1/(1-pz)}T_{|1-pz|}$ is symmetric. Thus, $(T_{1/(1-pz)}T_{|1-pz|})^{-1}=T_{1/|1-pz|}T_{p-z}=T_{\frac{p-z}{|1-pz|}}$ is symmetric with conjugation $C_{\psi_{p},\varphi_{p}}J$. \\ \par

In the following theorem, we find all complex symmetric weighted composition operators with conjugations $UJ$ that $U$ is unitary and $J$-symmetric weighted composition operator which was stated in Proposition 2.1.\\ \par

{\bf Theorem 2.5.} {\it Let $a_{0}\in \mathbb{D}$ and $a_{1},b\in \mathbb{C}$. Suppose that $\psi(z)=\frac{b}{1-a_{0}z}$ and $\varphi(z)=a_{0}+\frac{a_{1}z}{1-a_{0}z}$ that $\varphi$ is an analytic self-map of $\mathbb{D}$. \\
(1) For $p\neq 0$, the weighted composition operator  $C_{\widetilde{\psi},\widetilde{\varphi}}$ is complex symmetric with conjugation $C_{\psi_{p},\varphi_{p}}J$ if and only if $\widetilde{\psi}=\psi_{p}\cdot\psi \circ \varphi_{p}$ and $\widetilde{\varphi}=\varphi \circ \varphi_{p}$ for some $\varphi$ and $\psi$.\\
(2) For $|\lambda|=1$, the weighted composition operator  $C_{\widetilde{\psi},\widetilde{\varphi}}$ is complex symmetric with conjugation
$C_{\lambda z}J$ if and only if $\widetilde{\psi}=\psi(\lambda z)$ and $\widetilde{\varphi}(z)=\varphi(\lambda z)$ for some $\varphi$ and $\psi$.}\bigskip

{\bf Proof.}  (1) Let $\widetilde{\psi}=\psi_{p}\cdot \psi \circ \varphi_{p}$ and $\widetilde{\varphi}=\varphi \circ \varphi_{p}$ for some $\varphi$ and $\psi$. Then $C_{\widetilde{\psi},\widetilde{\varphi}}=T_{\psi_{p}\cdot \psi \circ \varphi_{p}}C_{\varphi \circ \varphi_{p}}=C_{\psi_{p},\varphi_{p}}C_{\psi,\varphi}$. Since $C_{\psi_{p},\varphi_{p}}$ and $C_{\psi,\varphi}$ are $J$-symmetric (see \cite[Proposition 2.9]{ham} and \cite[Theorem 3.3]{skl}), by Proposition 2.3, $C_{\widetilde{\psi},\widetilde{\varphi}}$ is $C_{\psi_{p},\varphi_{p}}J$-symmetric.
\par
Conversely, suppose that $C_{\widetilde{\psi},\widetilde{\varphi}}$ is complex symmetric with conjugation $C_{\psi_{p},\varphi_{p}}J$. By Proposition 2.3,  $C_{\psi_{p},\varphi_{p}}^{\ast}C_{\widetilde{\psi},\widetilde{\varphi}}$ is $J$-symmetric. The Cowen adjoint formula shows that $C_{\psi_{p},\varphi_{p}}^{\ast}$  is also a weighted composition operator. Then  $C_{\psi_{p},\varphi_{p}}^{\ast}C_{\widetilde{\psi},\widetilde{\varphi}}$ is a weighted composition operator which was defined in \cite[Proposition 2.9]{ham} and \cite[Theorem 3.3]{skl}. Then there is a $J$-symmetric weighted composition operator $C_{\psi,\varphi}$ that $C_{\widetilde{\psi},\widetilde{\varphi}}=C_{\psi_{p},\varphi_{p}}C_{\psi,\varphi}$. \par

(2) By the same idea which was stated in the proof of Part (1), the conclusion follows.\hfill $\Box$ \\
\par

{\bf Proposition 2.6.} {\it Suppose that $T$ is a bounded operator on a Hilbert space $H$. Then $T$ is a complex symmetric operator and an isometry if and only if $T$ is unitary.}\bigskip

{\bf Proof.}  Suppose that $T$ is complex symmetric with conjugation $C$. Since $T$ is an isometry, for each $f \in H^{2}$,
$$\|T^{\ast}Cf\|=\|CT^{\ast}Cf\|=\|Tf\|=\|f\|.$$
Then $\|T^{\ast}f\|=\|T^{\ast}C(Cf)\|=\|Cf\|=\|f\|$, and so $T^{\ast}$ is an isometry. We infer that $T$ is a unitary operator from
\cite[Proposition 2.17, p. 35]{c1} and \cite[Proposition 2.18, p. 35]{c1}.  \par
Conversely, it is obvious.\hfill $\Box$ \\ \par

In Theorem 2.7, we show that a weighted composition operator which is both a complex symmetric operator and an isometry is unitary; moreover, we find all conjugations for  unitary weighted composition operators.\\ \par

{\bf  Theorem 2.7.} {\it A weighted composition operator $C_{\psi,\varphi}$ is both an isometry and a complex symmetric operator if and only if $\varphi(z)=\lambda \frac{p-z}{1-\overline{p}z}$ and $\psi \equiv \psi_{p}$, where  $|\lambda|=1$ and $p \in \mathbb{D}$. Furthermore, if $p\neq 0$, then the conjugation for $C_{\psi,\varphi}$ is $C_{\psi_{p},\varphi_{p}}J$.}\bigskip

{\bf Proof.} By Proposition 2.6 and  \cite[Theorem 6]{Bourdon}, the first part is obvious.
According to \cite[Theorem 3.3]{skl} and \cite[Proposition 2.9]{ham}, $C_{\gamma z}$ is $J$-symmetric, when $|\gamma|=1$. Let $\gamma=\lambda\frac{p}{\overline{p}}$. We obtain $C_{\psi,\varphi}=C_{\psi_{p},\varphi_{p}}C_{\gamma z}$. By Propositions 2.1 and 2.3, we complete the proof. \hfill $\Box$ \\
\par

\section{Spectral theory}
 Recall that a nontrivial automorphism $\varphi$ of $\mathbb{D}$ (i.e., $\varphi$ is not the identity function of $\mathbb{D}$) is called
 elliptic if $\varphi$ has a fixed point in $\mathbb{D}$ and the other fixed point is in the complement of the closed disk.\par
 We say that $\varphi$ has a finite angular derivative at $\zeta \in \partial \mathbb{D}$ if the nontangential limit $\varphi(\zeta)$ exists, has modulus 1, and $\varphi'(\zeta)=\angle \lim_{z\rightarrow \zeta}\frac{\varphi(z)-\varphi(\zeta)}{z-\zeta}$ exists and finite.
Let $\varphi_{0}=I$ and $\varphi_{n}=\varphi\circ \varphi \circ ... \circ \varphi$ denote the $n$-th iterate of $\varphi$. If $\varphi$, not the identity and not an elliptic automorphism of $\mathbb{D}$, is a holomorphic self-map of $\mathbb{D}$, then there is a unique point $w$ in $\overline{\mathbb{D}}$ so that the iterates $\varphi_{n}$ of $\varphi$ tend to $w$ uniformly on compact subsets of $\mathbb{D}$ (see \cite[Theorem 2.51]{cm1}). The point $w$ will be referred to as the Denjoy-Wolff point of $\varphi$. We know that the Denjoy-Wolff point of $\varphi$  can be described as the unique fixed point of $\varphi$ in $\overline{\mathbb{D}}$ with $|\varphi'(w)|\leq 1$.\par
Suppose that $\varphi$, not an automorphism, is a linear-fractional self-map of $\mathbb{D}$ with a fixed point on $\partial \mathbb{D}$. Then $\varphi$ satisfies one of the following \par
(a) $\varphi$ is hyperbolic with one fixed point $\zeta \in \partial \mathbb{D}$ and the other fixed point outside the closed unit disk. Let $T(z)=\frac{\zeta+z}{\zeta-z}$. Then $\phi(z)=(T\circ\varphi\circ T^{-1})(z)=rz+t$, where $r=1/\varphi'(\zeta)$ (note that $r>1$) and $\mbox{Re}(t)\geq 0$ (and $\mbox{Re}(t)= 0$ if and only if $\varphi$ is an automorphism; moreover, in this case both fixed points of $\varphi$ lie on $\partial \mathbb{D}$). We call $t$ the translation number of $\varphi$. Then we obtain
\begin{equation}
\varphi(z)=\frac{(1+r-t)z+(r+t-1)\zeta}{(r-t-1)\overline{\zeta}z+(1+r+t)}.
\end{equation}
\par
(b) $\varphi$ is hyperbolic with one fixed point, $w$, inside the unit disk, and the other fixed point $\zeta$, on the unit circle. It is not hard to see that $\varphi$ is hyperbolic with this type (Denjoy-Wolff point of $\varphi$ is in $\mathbb{D}$) if and only if the Cowen auxiliary function $\sigma_{\varphi}$ is hyperbolic under the condition (a). Hence in this case
\begin{equation}
\varphi(z)=\frac{(1+r-\overline{t})z-(r-\overline{t}-1)\zeta}{1+r+\overline{t}-(r+\overline{t}-1)\overline{\zeta}z},
\end{equation}
where $t$ is the translation number of $\sigma_{\varphi}$ and $r=\frac{1}{\sigma'_{\varphi}(\zeta)}$. Note that in this case, since $\varphi$
has a Denjoy-Wolff point in $\mathbb{D}$, $\varphi$ is not an automorphism. Hence $\sigma_{\varphi}$ is not automorphism, so $\mbox{Re}(t)>0$.\par
(c) $\varphi$ is parabolic with only one fixed point $\zeta \in \partial \mathbb{D}$. Let $T(z)=\frac{\zeta+z}{\zeta-z}$. Then $\phi(z)=(T\circ\varphi\circ T^{-1})(z)=z+t$, where $\mbox{Re}(t)\geq 0$. Let us call $t$ the translation number of $\varphi$. Note $\mbox{Re}(t)= 0$
if and only if $\varphi$ is an automorphism. In \cite[p. 3]{sh} Shapiro showed that among the linear-fractional self-map of $\mathbb{D}$ fixing $\zeta \in
\partial \mathbb{D}$, the parabolic ones are characterized by
$\varphi'(\zeta)=1$.  We see that in this case
\begin{equation}
\varphi(z)=\frac{(2-t)z+t\zeta}{2+t-t\overline{\zeta}z}.
\end{equation}
Suppose that $\varphi_{1}$ and $\varphi_{2}$ are parabolic  with the same fixed point. It is not hard to see that $\varphi_{1} \circ \varphi_{2}$ is also parabolic. We use this fact in the proof of Theorem 3.5.\\ \par

{\bf  Lemma 3.1.} {\it Suppose that $\varphi$ is hyperbolic with fixed point $\zeta \in \partial \mathbb{D}$. If $\varphi$ is written as
\begin{equation}
\varphi(z)=a_{0}+\frac{a_{1}z}{1-a_{0}z},
\end{equation}
when $a_{0} \in \mathbb{D}$ and $a_{1}\in \mathbb{C}$, then $\varphi$ is an automorphism.}\bigskip

{\bf Proof.} First suppose that $\varphi$ is hyperbolic with Denjoy-Wolff point $\zeta \in \partial \mathbb{D}$. Then by Equation (3), we have
$$\varphi(z)=\frac{\frac{1+r-t}{1+r+t}z+\frac{(r+t-1)\zeta}{1+r+t}}{\frac{(r-t-1)\overline{\zeta}}{1+r+t}z+1}.$$

If $\varphi$ is written as in  Equation (6), then $\frac{-r+t+1}{1+r+t}\overline{\zeta}=\frac{r+t-1}{1+r+t}\zeta$. Therefore, $-r\overline{\zeta}+t\overline{\zeta}+\overline{\zeta}=r\zeta+t\zeta-\zeta$ and so $(\mbox{Re}(\zeta))(1-r)=t(\mbox{Im} (\zeta))i$.
It shows that $t$ is pure imaginary and the result follows.  Now suppose that $\varphi$ is hyperbolic with Denjoy-Wolff $w \in \mathbb{D}$ and fixed point $\zeta \in \partial \mathbb{D}$. By Equation (4), assume that
$$\varphi(z)=\frac{\frac{(1+r-\overline{t})}{1+r+\overline{t}}z-\frac{r-\overline{t}-1}{1+r+\overline{t}}\zeta}{1-\frac{(r+\overline{t}-1)\overline{\zeta}z}{1+r+\overline{t}}}.$$
 Since $\varphi$ is as in Equation (6), $\frac{(r+\overline{t}-1)\overline{\zeta}}{1+r+\overline{t}}=-\frac{(r-\overline{t}-1)\zeta}{1+r+\overline{t}}$. Then $(\mbox{Re}(\zeta))(r-1)=\overline{t}(\mbox{Im} (\zeta))i$.
 It follows that $t$ is pure imaginary which is a contradiction.\hfill $\Box$ \\

{\bf  Lemma 3.2.} {\it Suppose that $\varphi$ is parabolic with fixed point $\zeta \in \partial \mathbb{D}$. Then $\varphi$ is as in Equation (6) if and only if $\zeta=1$ or $\zeta=-1$.}\bigskip

{\bf Proof.} Suppose that $\varphi$ is written as in Equation (6) and by Equation (5), $\varphi(z)=\frac{\frac{(2-t)z}{2+t}+\frac{t\zeta}{2+t}}{-\frac{t\overline{\zeta}}{2+t}z+1}$. If $t=0$, then $\varphi(z)=z$ and $\varphi$ is not parabolic. Then we assume that  $t\neq 0$. Since $\varphi$ is as in Equation (6), $\overline{\zeta}=\zeta$. It shows that $\zeta=1$ or $\zeta=-1$. \par
Conversely, it is obvious. \hfill $\Box$ \\

Note that by the proof of the pervious lemma, we see that if $\varphi$   is parabolic with Denjoy-Wolff point $1$ which is written as in Equation (5), then
\begin{equation}
\varphi(z)=\frac{\frac{2-t}{2+t}z+\frac{t}{2+t}}{\frac{-t}{t+2}z+1}.
\end{equation}

{\bf Lemma 3.3.} {\it Suppose that $p\in  \mathbb{D}-\{0\}$. If $1$ is the fixed point of $\varphi_{p}$, then $\varphi_{p}$ is a parabolic automorphism. Moreover, if $-1$ is the fixed point of $\varphi_{p}$, then $\varphi_{p}$ is a hyperbolic automorphism with Denjoy-Wolff point $-1$.}\bigskip

{\bf Proof.} Assume that $\varphi_{p}(1)=1$. Note that  $p=|p|\cos(\theta)+|p|\sin(\theta)i$, where $\theta=\mbox{Arg}(p)$. Then $\overline{p}(p-1)=p(1-\overline{p})$ and so $|p|^{2}=\mbox{Re}(p)=|p|\cos(\theta)$. Hence $|p|=\cos(\theta)$. Since $\frac{\overline{p}}{p}=e^{-2\theta i}$, \cite[Exercise 4, p. 7]{sh} implies that $\varphi_{p}$ is parabolic. Now suppose that $\varphi_{p}(-1)=-1$. Then $\overline{p}(p+1)=-p(1+\overline{p})$ and so $|p|^{2}=-\mbox{Re}(p)=-|p|\cos(\theta)$. Hence $|p|=-\cos(\theta)$. Again by  \cite[Exercise 4, p. 7]{sh}, $\varphi_{p}$ is a hyperbolic automorphism. Now we show that $-1$ is the Denjoy-Wolff point of $\varphi_{p}$. We have $\varphi'_{p}(-1)=\frac{\overline{p}}{p}\frac{|p|^{2}-1}{(1+\overline{p})^{2}}$. We know that $p=|p|\cos(\theta)+|p|\sin(\theta)i=-\cos^{2}(\theta)-\sin(\theta)\cos(\theta)i$. Therefore, $1+\overline{p}=1-\cos^{2}(\theta)+\sin(\theta)\cos(\theta)i=\sin^{2}(\theta)+\sin(\theta)\cos(\theta)i$.  We obtain $|\varphi'_{p}(-1)|^{2}=|\frac{|p|^{2}-1}{(1+\overline{p})^{2}}|^{2}=|\frac{\cos^{2}(\theta)-1}{(\sin^{2}(\theta)+\sin(\theta)\cos(\theta)i)^{2}}|^{2}=
\frac{\sin^{4}(\theta)}{\sin^{4}(\theta)+
\sin^{2}(\theta) \cos^{2}(\theta)} < 1$. Then $|\varphi'_{p}(-1)|<1$ and it follows that $-1$ is the Denjoy-Wolff point of $\varphi_{p}$.  \hfill $\Box$ \\

In the rest of this paper, we suppose that $C_{\psi,\varphi}$ is $J$-symmetric and    $\varphi$ and $\psi$ were represented in Theorem 2.5. From now on, unless otherwise stated, we assume that $C_{\widetilde{\psi},\widetilde{\varphi}}$ is weighted composition operator which was given in the first part of Theorem 2.5. Suppose that $T$ is a bounded operator on a Hilbert space $H$.  Through this paper, the spectrum of $T$ and the spectral
radius of $T$ are denoted by $\sigma(T)$ and $r(T)$, respectively.\\ \par

{\bf Theorem  3.4.} {\it If $C_{\widetilde{\psi},\widetilde{\varphi}}$ is compact or power compact, then $r(C_{\widetilde{\psi},\widetilde{\varphi}})=|\psi_{p}(w)\psi(\varphi_{p}(w))|$ and $\sigma(C_{\widetilde{\psi},\widetilde{\varphi}})=\{\psi_{p}(w)\psi(\varphi_{p}(w))(\widetilde{\varphi}'(w))^{m}:m=0,1,...\}\cup\{0\}$, where $w$ is the Denjoy-Wolff point of $\widetilde{\varphi}=\varphi\circ\varphi_{p}$.}\bigskip

{\bf Proof.} If $\varphi \circ\varphi_{p}$ is compact or power compact, then it is easy to see that $\varphi \circ \varphi_{p}$ has a Denjoy-Wolff point $w \in  \mathbb{D}$. There is an integer $n$ such that $C_{\widetilde{\psi},\widetilde{\varphi}}^{n}=C_{\widetilde{\psi}\cdot\widetilde{\psi}\circ\widetilde{\varphi}...\widetilde{\psi}\circ\widetilde{\varphi}_{n-1},\widetilde{\varphi}_{n}}$ is compact. By the Spectral Mapping Theorem, $\sigma_{e}(C_{\widetilde{\psi},\widetilde{\varphi}})=\{0\}$. Moreover, by \cite[Theorem 1]{g1}, $\sigma(C_{\widetilde{\psi},\widetilde{\varphi}}^{n})=\{(\widetilde{\psi}(w))^{n}(\widetilde{\varphi}'(w))^{mn}:m=0,1,...\}$ and it follows from the Spectral Mapping Theorem that all elements of $\sigma(C_{\widetilde{\psi},\widetilde{\varphi}})$ are in $\partial \sigma(C_{\widetilde{\psi},\widetilde{\varphi}})$. Then by \cite[Proposition 6.7, p. 210]{c1} and \cite[Proposition 4.4, p. 359]{c1},
 $\sigma(C_{\widetilde{\psi},\widetilde{\varphi}})=\sigma_{e}(C_{\widetilde{\psi},\widetilde{\varphi}}) \cup \sigma_{p}(C_{\widetilde{\psi},\widetilde{\varphi}})$. The result follows by \cite[Proposition 2.6]{ham} and \cite[Theorem 4.3]{skl}.\hfill $\Box$ \\

Now suppose that $\widetilde{\varphi}=\varphi\circ\varphi_{p}$ is not power compact and it is not an automorphism. Then either\\
(i) $\varphi(\zeta)=\zeta $ and $\varphi_{p}(\zeta)=\zeta $, where $\zeta \in \partial \mathbb{D}$\\
or\\
(ii) $\varphi(\zeta)=\eta $ and $\varphi_{p}(\eta)=\zeta $, where $\zeta,\eta \in \partial \mathbb{D}$ and $\zeta \neq \eta$.\par
By Lemmas 3.1 and 3.2, we see that if $\varphi$ and $\varphi_{p}$ satisfy the conditions of Part (i), then either $\zeta=1$ or $\zeta=-1$. In the following  theorem, we find $r(C_{\widetilde{\psi},\widetilde{\varphi}})$, when $\varphi$ and $\varphi_{p}$  satisfy the conditions of Part (i) and $\widetilde{\varphi}$ is parabolic.\\ \par

{\bf Theorem  3.5.} {\it Suppose that $\widetilde{\varphi}$ is not an automorphism. If $\varphi$ and $\varphi_{P}$  fix $1$, then  $\varphi,\varphi_{p}$ and $\widetilde{\varphi}$ are parabolic, $r(C_{\widetilde{\psi},\widetilde{\varphi}})=|\frac{\psi(0)(1-|p|^{2})^{1/2}(2+t)}{2(1-\overline{p})}|$ and $\sigma(C_{\widetilde{\psi},\widetilde{\varphi}})=\{\frac{\psi(0)(1-|p|^{2})^{1/2}(2+t)}{2(1-\overline{p})} e^{-b(t+\widetilde{t})}: b\geq 0\} \cup \{0\}$, where $t$ and $\widetilde{t}$ are  the translation number of $\varphi$ and $\varphi_{p}$, respectively.}\bigskip

{\bf Proof.} Since $\varphi_{p}$ is an automorphism and $\varphi\circ\varphi_{p}$ is not an automorphism, $\varphi$ is not an automorphism. Invoking Lemmas 3.1 and 3.2,  $\varphi$ must be parabolic with fixed point $1$.
  By Lemma 3.3, $\varphi_{p}$ must be parabolic with fixed point $1$. Hence $\widetilde{\varphi}=\varphi\circ\varphi_{p}$ is also parabolic.  Since $\varphi$ is parabolic with fixed point $1$ and $C_{\psi,\varphi}$ is $J$-symmetric, by Equation (7), we find that $\psi(z)=\frac{\psi(0)}{1-\frac{t}{2+t}z}$. Then by \cite[Theorem 4.7]{Bourdon1}, $r(C_{\widetilde{\psi},\widetilde{\varphi}})=|\widetilde{\psi}(1)|= |\psi_{p}(1)\cdot \psi(1)|=|\frac{\psi(0)(1-|p|^{2})^{1/2}(2+t)}{2(1-\overline{p})}|$ and  $\sigma(C_{\widetilde{\psi},\widetilde{\varphi}})=\{\frac{\psi(0)(1-|p|^{2})^{1/2}(2+t)}{2(1-\overline{p})} e^{-b(t+\widetilde{t})}: b\geq 0\} \cup \{0\}$.\hfill $\Box$ \\

In the next theorem, we find $r(C_{\widetilde{\psi},\widetilde{\varphi}})$, when $\widetilde{\varphi}$ is a hyperbolic non-automorphism with Denjoy-Wolff point $-1$.\\ \par

{\bf Theorem  3.6.} {\it Suppose that $\widetilde{\varphi}$ is not an automorphism. If $\varphi$ and $\varphi_{p}$  fix  $-1$, then $\widetilde{\varphi}$ is hyperbolic with Denjoy-Wolff point $-1$, $\varphi$ is parabolic with translation number $t$, $r(C_{\widetilde{\psi},\widetilde{\varphi}})=|\frac{2+t}{2(1+\overline{p})}(\frac{1-|p|^{2}}{\varphi'_{p}(-1)})^{1/2}|$ and $\sigma(C_{\widetilde{\psi},\widetilde{\varphi}})=\{z:|z| \leq r(C_{\widetilde{\psi},\widetilde{\varphi}})\}$.}\bigskip

{\bf Proof.} It is easy to see that $-1$ is the fixed point of $\varphi \circ \varphi_{p}$. We see that $(\varphi\circ\varphi_{p})'(-1)=\varphi'(\varphi_{p}(-1))\cdot \varphi'_{p}(-1)$. Since by Lemmas 3.1 and 3.2, $\varphi$ is parabolic, $\varphi'(-1)=1$ and so $(\widetilde{\varphi})'(-1)=\varphi'_{p}(-1)$. Lemma 3.3 shows that $|(\varphi \circ\varphi_{p})'(-1)|< 1$ and $-1$ is the Denjoy-Wolff point of $\widetilde{\varphi}$. Then $\widetilde{\varphi}$ must be a hyperbolic non-automorphism with Denjoy-Wolff point $-1$. From Equation (5), $\varphi(z)=\frac{(2-t)z-t}{tz+2+t}=\frac{\frac{2-t}{2+t}z-\frac{t}{2+t}}{1+\frac{t}{2+t}z}$, where $t$ is the translation number of $\varphi$. Since $C_{\psi,\varphi}$ is $J$-symmetric, it is easy to see that $\psi(z)=\frac{1}{1+\frac{t}{2+t}z}$. Invoking \cite[Theorem 4.5]{Bourdon1}, $r(C_{\widetilde{\psi},\widetilde{\varphi}})=|\frac{\psi(-1)}{1+\overline{p}}(\frac{1-|p|^{2}}{\varphi'_{p}(-1)})^{1/2}|=|\frac{2+t}{2(1+\overline{p})}(\frac{1-|p|^{2}}{\varphi'_{p}(-1)})^{1/2}|$ and $\sigma(C_{\widetilde{\psi},\widetilde{\varphi}})=\{z:|z| \leq r(C_{\widetilde{\psi},\widetilde{\varphi}})\}$.\hfill $\Box$ \\

Note that although every complex symmetric composition $C_{\varphi}$ must have a fixed point in $\mathbb{D}$ (see \cite[Proposition 2.1]{bnoor}), Theorems 3.5 and 3.6 showed that there are complex symmetric weighted composition operators $C_{\widetilde{\psi},\widetilde{\varphi}}$ that $\widetilde{\varphi}$ has no fixed point in $\mathbb{D}$. \par
In \cite[Theorem 4.6]{skl} Jung et al. obtained an inequality for $r(C_{\psi,\varphi})$, when $C_{\psi,\varphi}$ is $J$-symmetric.
In the next corollary, we find the spectral radius of $J$-symmetric weighted composition operators, when $\varphi$ is not an automorphism. Note that for an automorphism $\varphi$, the spectral radius of $C_{\psi,\varphi}$ was found in \cite{hlns}.\\ \par

{\bf Corollary  3.7.} {\it  Assume that $\varphi$ is not an automorphism. Let $w$ be the Denjoy-Wolff point of $\varphi$.\\
(a) If $w \in \mathbb{D}$, then $r(C_{\psi,\varphi})=|\psi(w)|$.\\
(b) If $w \in \partial \mathbb{D}$, then $\varphi$ is parabolic and $w$ is either $1$ or $-1$. Moreover, if $w=1$, then $r(C_{\psi,\varphi})=\left|\frac{\psi(0)(2+t)}{2}\right|$  and if $w=-1$, then $r(C_{\psi,\varphi})=\left| \frac{\psi(0)(2+t)}{2+2t}\right|$, where $t$ is the translation number of $\varphi$.}\bigskip

{\bf Proof.} (a) If $w \in \mathbb{D}$, then by Lemmas 3.1 and 3.2, $\varphi$ has no fixed point in $\partial \mathbb{D}$. Then $C_{\varphi}$ is compact or power compact. The result follows by the Spectral Mapping Theorem, \cite[Theorem 4.3]{skl} and \cite[Proposition 2.6]{ham}.\\
(b) Assume that $\varphi$ has a Denjoy-Wolff point $w \in \partial \mathbb{D}$. As we saw in Lemmas 3.1 and 3.2, $\varphi$ must be parabolic with fixed point $1$ or $-1$. By \cite[Theorem 4.7]{Bourdon1} and the similar idea which was stated in the proof of Theorem 3.5, we see that if $w$=1, then $r(C_{\psi,\varphi})=|\psi(1)|=\left|\frac{\psi(0)(2+t)}{2}\right|$ and if $w=-1$, then $r(C_{\psi,\varphi})=|\psi(-1)|=\left| \frac{\psi(0)(2+t)}{2+2t}\right|$.\hfill $\Box$ \\

It is not hard to see that for $|\lambda|=1$, $C_{\varphi(\lambda z)}$ is not power compact if and only if there is $\zeta \in \partial \mathbb{D}$ such that $\varphi(\lambda \zeta)=\zeta$. In the following proposition, we  find the spectrum of power compact weighted composition operator $C_{\widetilde{\psi},\widetilde{\varphi}}$, when $C_{\widetilde{\psi},\widetilde{\varphi}}$ was given in the second part of Theorem 2.5.\\ \par

{\bf Proposition  3.8.} {\it Assume that $\varphi$ is not an automorphism. Suppose that $C_{\widetilde{\psi},\widetilde{\varphi}}$ is as in the second part of Theorem 2.5. Let $\lambda$ be constant and $|\lambda|=1$. If for each $\zeta \in \partial \mathbb{D}$, $\varphi(\lambda \zeta) \neq \zeta$, then $\sigma(C_{\widetilde{\psi},\widetilde{\varphi}})=\{\psi(\lambda w)(\lambda \varphi'(\lambda w))^{n}: n=0,1,...\}$ and $r(C_{\widetilde{\psi},\widetilde{\varphi}})=|\psi(\lambda w)|$, when $w \in \mathbb{D}$ and $\varphi(\lambda w)=w$.}\bigskip

{\bf Proof.} One may easily see that in this case $C_{\widetilde{\psi},\widetilde{\varphi}}$ is power compact and the result follows by the similar idea which was stated in the proof of Theorem 3.4.\hfill $\Box$ \\

\bigskip
{M. Fatehi, Department of Mathematics, Shiraz Branch, Islamic Azad
University, Shiraz, Iran. \par E-mail: fatehimahsa@yahoo.com \par
}

\end{document}